\documentclass[12pt]{article}
\usepackage{amsmath,amssymb,pb-diagram}%,pdfsync}

\def \1{{\bf 1}}
\def \a{{{\mathfrak a}}}

\def \Ad{{\rm Ad}}

\def \bs{\backslash}

\def \C{{\mathbb C}}

\def \CC{{\cal C}}

\def \CE{{\cal E}}

\def \CH{{\cal H}}

\def \CO{{\cal O}}

\def \det{{\rm det}}
\def \df{\ \begin{array}{c} _{\rm def}\\ ^{\displaystyle =}\end{array}\ }
\def \diag{{\rm diag}}
\def \ds{\displaystyle}

\def \End{{\rm End}}

\def \eqn{\begin{eqnarray*}}
\def \endeqn{\end{eqnarray*}}
\def \F{{\mathbb F}}

\def \ga{\gamma}
\def \Ga{\Gamma}

\def \GL{{\rm GL}}

\def \La{\Lambda}
\def \m{{{\mathfrak m}}}
\def \mod{{\rm mod}}

\def \N{\mathbb N}

\def \ph{\varphi}
\def \prf{\noindent{\bf Proof: }}

\def \PGL{{\rm PGL}}

\def \qed{\ifmmode\eqno $\square$\else\noproof\vskip 12pt plus 3pt minus 9pt \fi}
 \def\noproof{{\unskip\nobreak\hfill\penalty50\hskip2em\hbox{}%
     \nobreak\hfill $\square$\parfillskip=0pt%
     \finalhyphendemerits=0\par}}
\def \Q{\mathbb Q}

\def \SL{{\rm SL}}

\def \tr{{\hspace{1pt}\rm tr\hspace{2pt}}}
\def \vol{{\rm vol}}

\def \Z{\mathbb Z}
\def \={\ =\ }
\def \({\left(}
\def \){\right)}

\newtheorem{theorem}{Theorem}[section]

\newtheorem{lemma}[theorem]{Lemma}

\newtheorem{proposition}[theorem]{Proposition}

\renewcommand{\sp}[1]{\left\langle #1 \right\rangle}

\begin{document}

\pagestyle{myheadings} \markright{IHARA-SELBERG...}

\title{The Ihara-Selberg zeta function for $\PGL_3$ and Hecke operators}
\author{Anton Deitmar \& J. William Hoffman}

\date{}
\maketitle

{\bf Abstract.}
A weak version of the Ihara formula is proved for zeta functions attached to quotients of the Bruhat-Tits building of $\PGL_3$.
This formula expresses the zeta function in terms of Hecke-Operators. It is the first step towards an arithmetical interpretation of the combinatorially defined zeta function.

\tableofcontents

\newpage
\section*{Introduction}
Y. Ihara \cite{Ihara} extended the theory of Selberg type zeta functions to $p$-adic settings.
His work was later generalized by K. Hashimoto \cite{Hash1,Hash2,Hash3}, H. Bass \cite{Bass},
H. Stark and A. Terras \cite{ST}, and others.
Ihara defined the zeta function in group theoretical terms first, but it can be described geometrically as follows.
Let $\Ga\bs X$ be a finite quotient of the Bruhat-Tits tree of a rank one $p$-adic group modulo an arithmetic group $\Gamma$.
Then define the zeta function by
$$
Z(u)\=\prod_c (1-u^{l(c)}),
$$ 
where the product runs over all primitive closed loops in $\Ga\bs X$.
Ihara proved the remarkable formula
$$
Z(u)\= (1-u^2)^{-\chi} \det( 1-Au + qu^2),
$$
where $A$ is the adjacency operator on $\Ga\bs X$ which can be interpreted as the canonical generator of the unramified Hecke-algebra.
Further, $\chi<0$ is the Euler-characteristic of $\Ga\bs X$, and $q$ is the order of the residue class field.

For $\Ga$ being the unit group of the maximal order in a quaternion algebra, this formula allowed Ihara to relate $Z(u)$ to the Hasse-Weil zeta function of the Shimura curve attached to $\Gamma$.
This is the only proven link between Selberg-type zeta functions and arithmetical zeta functions.

In \cite{padgeom} the author gave a definition of an Ihara-type zeta function $Z(u)$ for a higher rank group.
There is no Ihara-formula for higher rank up to date.
In this paper we give an approximation to an Ihara formula in the case of the group $\PGL_3$.
For this group the unramified Hecke-algebra has two generators $\pi_1,\pi_2$.
The canonical replacement of the determinant factor in Ihara's formula is
$$
\det(1-u\pi_1 +u^2 q \pi_2 -u^3 q^3).
$$
The main result of the present paper is

\begin{theorem}\label{main}
There are a natural number $n$ and a polynomial $P(u)$ such that
$$
Z(u)\= \frac{\det(1-u\pi_1+ u^2 q \pi_2 -u^3 q^3)^n}{P(u)}.
$$             
\end{theorem}

I thank A. Setyadi for pointing out an error in an earlier version.

\section{The building}
Let $F$ be a non-archimedean local field.
Let $\CO$ be its valuation ring with maximal ideal $\m\subset\CO$.
Fix a generator $\varpi$ of $\m$ and let $q$ be the cardinality of the residue class field $k=\CO/\m$.

Consider the locally compact group $G=\PGL_3(F)= \GL_3(F)/F^\times$.
It is totally disconnected and every maximal compact subgroup is conjugate to $K=\PGL_3(\CO)=\GL_3(\CO)/\CO^\times$.
Let $X$ be the Bruhat-Tits building \cite{Brown} of $G$.
In this particular case the Bruhat-Tits building can be described rather nicely.
The vertex set $X_0$ of $X$ is the set of homothety classes of $\CO$-lattices in $F^3$.
Recall that an $\CO$-lattice in $F^3$ is a finitely generated $\CO$-submodule $\La$ of $F^3$ such that $F\La=F^3$.
Two lattices $\Lambda, \Lambda'$ are \emph{homothetic}, if there exists $\alpha\in F^\times$ such that $\Lambda'=\alpha\Lambda$.
Every lattice $\La$ is the image under some $g\in\GL_3(F)$ of the standard lattice $L_0=\CO e_1\oplus \CO e_2\oplus\CO e_3$, where $e_1,e_2,e_3$ is the standard basis of $F^3$.
The set of all lattices thus can be identified with $\GL_3(F)/\GL_3(\CO)$ and the set $X_0$ of homothety classes of lattices with $G/K$.
Let $G'$ denote the image of $\SL_3(F)$ in $G$.
The set $X_0$ splits into three orbits under the action of $G'$.
These orbits are given by $L_0$ as above, $L_1=\CO e_1\oplus\CO e_2\oplus\varpi\CO e_3$, and $L_2=\CO e_1\oplus\varpi\CO e_2\oplus\varpi\CO e_3$.
For a given vertex $x\in X_0$ we say $x$ \emph{is of type} $j$ if $G' x=G' L_j$ for $j=0,1,2$.
Two vertices $x\ne y$ are joined by an edge if and only if there are representatives $\La_1$ and $\La_2$ for $x$ and $y$ such that $\varpi\La_1\subset\La_2\subset\La_1$.
It follows that $x$ and $y$ must be of different type.
This describes the $1$-skeleton $X_1$ of $X$.
The following Lemma gives further properties of the graph $X_1$.

\begin{lemma}\label{graph}
\begin{enumerate}
\item 
Every vertex in $X$ has $2(q^2+q+1)$ neighbours.
\item
Two neighboured vertices have $q+1$ common neighbours.
\item
Any three distinct vertices have at most one common neighbour.
\end{enumerate}
\end{lemma}

\prf
For (a) it suffices to consider the vertex given by $L_0$.
Every neighbour has a representative lattice $L$ with
$$
\varpi L_0\ \subset\ L\ \subset\ L_0.
$$
Now $L_0/\varpi L_0\cong \F_q^3$ as a vector space over $\F_q$, and $L$ defines a sub vector space.
Thus the set of all neighbours of $[L_0]$ is in bijection with the set of all non-trivial sub vector spaces of $\F_q^3$ which are $2(1+q+q^2)$ in number.
Part (b) and (c) are similar.
\qed

Whenever three vertices $x,y,z$ are mutually connected by edges, then this triangle forms the boundary of a 2-cell of $X$, called a \emph{chamber}.
This describes $X$ as a CW-complex.
There is, however, more structure through the geometry of the apartments.
For instance, whenever two edges meet in a vertex, there is an angle between them which can be $\pi/3, 2\pi/3,$ or $\pi$.
A \emph{geodesic} $c$ in $X$ is a straight oriented line in one apartment.
If $c$ happens to lie inside $X_1$, then it gives rise to a sequence of edges $(\dots,e_{-1},e_0,e_1,\dots)$ such that $e_k$ and $e_{k+1}$ have angle $\pi$ for every $k\in\Z$.
In this case we say that $c$ is a \emph{rank-one geodesic}.

For each edge $e$ with vertices $\{x,y\}$ we fix an orientation, i.e., an ordering of the vertices $(x,y)$ such that if $x$ is of type $j$, then $y$ is of type $j+1\ \mod(3)$.
An edge equipped with this orientation will be called \emph{positively oriented}.
Likewise, for the chambers we fix a positive orientation by ordering the vertices by type.

\begin{lemma}\label{orientation}
The action of $G$ on the edges and the chambers preserves the positive orientation.
\end{lemma}

\prf
Let $g\in G$.
It suffices to show that if $gL_0$ is of type $j$, then $gL_1$ is of type $j+1\ \mod(3)$.
First note that the double quotient $G'\bs G/K$ has three elements given by the class of $1$, $\diag(1,1,\varpi)$ and $\diag(1,\varpi,\varpi)$.
Next note that the action of $K$ preserves the positive orientation on edges that contain the base point $L_0$.
Thus it suffices to prove the claim for the three given elements which is easily done.
\qed

\section{The zeta function}
Let $\Ga\subset G$ be a discrete cocompact and torsion-free subgroup.
Then $\Ga$ acts without fixed points on $X$ and thus $\Ga$ is the fundamental group of the quotient $\Ga\bs X$.

\subsection{Definition}
A \emph{geodesic} $c$ in the quotient $\Ga\bs X$ is the image of a geodesic $\tilde c$ in $X$ under the projection map $X\to\Ga\bs X$.
The geodesic $c$ is called \emph{rank-one} if $\tilde c$ is.
For a geodesic $c$ we denote by $c^{-1}$ the geodesic with the reversed orientation.
When speaking about \emph{closed geodesics} in $\Ga\bs X$ we adopt the convention that a closed geodesic comes with a multiplicity (going round more then once).
A closed geodesic with multiplicity one is called a \emph{primitive} closed geodesic.
To a given closed geodesic $c$ there is a unique primitive one $c_0$ such that $c$ is a power of $c_0$.
For a closed geodesic $c$ in $\Ga\bs X$ let $l(c)$ denote its length.
Here the length is normalized such that any edge gets the length 1.
We define the \emph{zeta function}
$$
Z(u)\=\prod_c\left( 1-u^{l(c)} \right)
$$
as a formal power series at first, where $c$ ranges over the set of all primitive rank-one closed geodesics in $\Ga\bs X$ modulo homotopy and modulo change of orientation.
It is easy to show that the Euler product defining $Z(u)$ actually converges for $u\in\C$ with $|u|$ small enough.

\subsection{A comparison}
An element $g$ of $G$ is called \emph{neat} if for every rational representation $\rho\colon G\to \GL_n(F)$ over $F$ the matrix $\rho(g)$ has the following property: the subgroup of $\bar F^\times$ generated by all eigenvalues of $\rho(g)$ is torsion-free. 
Here $\bar F$ is an algebraic closure of $F$.
The element $g$ is called \emph{weakly neat} if the adjoint $\Ad(g)\in\GL({\rm Lie}(G))$ has no non-trivial root of unity as eigenvalue.
Obviously neat implies weakly neat.
A subgroup $\Ga\subset G$ is called \emph{neat}/\emph{weakly neat} if every $\ga\in\Ga$ is neat/weakly neat in $G$.
Every arithmetic group $\Ga$ has a subgroup of finite index which is neat \cite{Borel}.

An element $g$ of $G$ is called \emph{regular} if its centralizer is a torus.
A subgroup $\Ga$ of $G$ is called regular if every $\ga\in\Ga$, $\ga\ne 1$ is regular in $G$.
A regular group is weakly neat.

In \cite{padgeom}, the author defined for $\Ga$ being discrete, cocompact, and weakly neat a zeta function $Z_P(u)$ attached to  a parabolic subgroup $P\subset G$ of splitrank one.
It is shown that $Z_P(u)$ is a rational function and that its poles and zeros can be described in terms of certian cohomology groups.

\begin{proposition}
Suppose the group $\Ga$ is discrete, cocompact and regular.
Then $Z_P(u)=Z(u)$ for every parabolic $P$ of splitrank one.
\end{proposition}

\prf
We will recall the definition of $Z_P(u)$.
Let $P=LN$ be a Levi decomposition of P and $A\subset L$ be a maximal split torus.
The dimension of $A$ is one.
Let $A^+$ be the set of all $a\in A$ that act on the Lie algebra of $N$ by eigenvalues $\mu$ with $|\mu|>1$.
Fix an isomorphism $\ph\colon A\cong F^\times$ that maps $A^+$ to the set of $x\in F^\times$ with $v(x)>0$, where $v$  is the valuation of $F$.
For $a\in A^+$ let $l(a)=v(\ph(a))$.
Let $M$ be the derived group of $M$ and let $M_{\rm ell}$ be the set of all elliptic elements of $M$.
Let $\CE_P(\Ga)$ denote the set of all conjugacy classes $[\ga]$ in $\Ga$ such that $\ga$ is in $G$ conjugate to an element $a_\ga m_\ga\in A^+ M_{\rm ell}$.
An element $\ga\in\Ga$ is called \emph{primitive} if $\ga=\sigma^n$m $n\in\N$, $\sigma\in\Ga$ implies $n=1$.
Let $\CE_P^p(\Ga)$ denote the set of primitive elements in $\CE_P(\Ga)$.
The zeta function $Z_P$ is defined as
$$
Z_P(u)\=\prod_{[\ga]\in\CE_P^p(\Ga)}(1-u^{l(a_\ga)})^{\chi_1(\Ga_\ga)},
$$
where  $\Ga_\ga$ is the centralizer of $\ga$ in $\Ga$ and
$$
\chi_1(\Ga_\ga)\=\sum_{p=0}^{\dim X} p(-1)^{p+1}\dim H^p(\Ga_\ga,\Q).
$$
First we remark that since $\Ga$ is regular, we have that $\Ga_\ga\cong\Z$ for every $\ga\in\CE_P(\Ga)$ and thus the Euler numbers $\chi_1(\Ga_\ga)$ are all equal to $1$.
Next let $[\ga]\in\CE_P^p(\Ga)$.
The function $d_\ga(x)={\rm dist}(\ga x,x)$ on $X$ attains its minimum on a unique apartment of $X$.
On this apartment, $\ga$ acts by a translation along a rank-one geodesic $\tilde c$ by the amount $l(a_\ga)$.
So $\ga$ closes this geodesic and its image $c$ in $\Ga\bs X$ has length $l(c)=l(a_\ga)$.
The other way round, every rank-one closed geodesic $c$ must be closed by one primitive element $\ga$ of $\Ga$.
Then $\ga$ either lies in $\CE_P(\Ga)$ or has splitrank two in which case the geodesic it closes cannot be rank-one.  
This shows that the Euler products defining $Z(u)$ and $Z_P(u)$ coincide.
\qed

\subsection{A factorization}
Two rank-one geodesics in $X$ are called \emph{adjacent} if they lie in the same apartment, they are parallel, and there is only one row of chambers between them.
Recall a \emph{gallery} \cite{Brown} in $X$ is a sequence $g=(C_0,\dots, C_n)$ of chambers such that $C_j$ and $C_{j+1}$ are adjacent for every $j$.
We say that a gallery $g$ is \emph{rank-one} if $C_{j-1}\ne C_{j+1}$ for every $j$ and the gallery is located between two rank-one geodesics.
The next picture shows an example of a rank-one gallery.

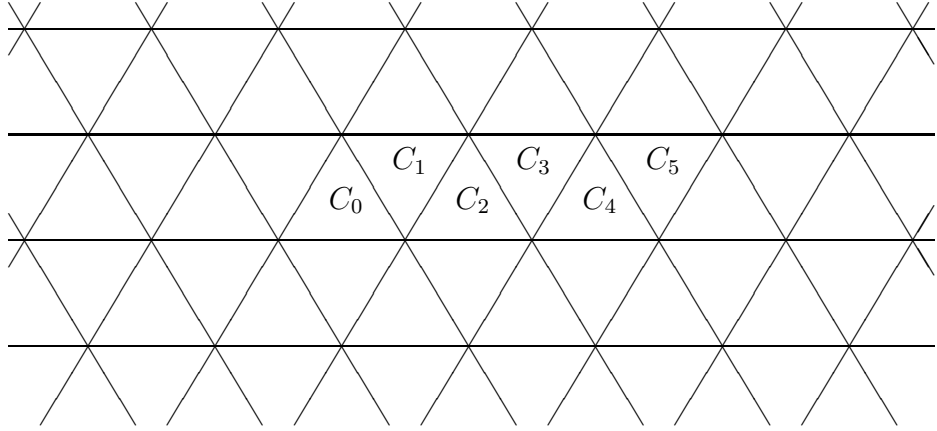
\begin{picture}(300,160)
% waagrecht
\put(0,30){\line(1,0){350}}
\put(0,70){\line(1,0){350}}
\put(0,110){\line(1,0){350}}
\put(0,150){\line(1,0){350}}
% fallend
\put(0,80){\line(3,-5){48}}
\put(0,160){\line(3,-5){96}}
\put(48,160){\line(3,-5){96}}
\put(96,160){\line(3,-5){96}}
\put(144,160){\line(3,-5){96}}
\put(192,160){\line(3,-5){96}}
\put(240,160){\line(3,-5){96}}
\put(288,160){\line(3,-5){62}}
\put(336,160){\line(3,-5){14}}
% steigend
\put(0,140){\line(3,5){12}}
\put(0,60){\line(3,5){60}}
\put(12,0){\line(3,5){96}}
\put(60,0){\line(3,5){96}}
\put(108,0){\line(3,5){96}}
\put(156,0){\line(3,5){96}}
\put(204,0){\line(3,5){96}}
\put(252,0){\line(3,5){96}}
\put(300,0){\line(3,5){50}}

% Gallerie
%\put(121,158){$C_0$}
\put(121,82){$C_0$}
\put(145,97){$C_1$}
\put(169,82){$C_2$}
\put(192,97){$C_3$}
\put(217,82){$C_4$}
\put(241,97){$C_5$}
% Ende
\end{picture}

A rank-one gallery in $\Ga\bs X$ is the image of a rank-one gallery in $X$ under the projection map.
In $\Ga\bs X$ it may happen for a rank-one gallery $g=(C_0,\dots,C_n)$ that $C_0=C_n$ in which case we say that $g$ is \emph{closed}.
In this case the number $n$ is even and we define the \emph{length} of $g$ to be $l(g)=n/2$.
We say that $g$ is \emph{primitive} if furthermore $C_0\ne C_j$ for $0<j<n$.
Two closed galleries $(C_0,\dots,C_n)$ and $(E_0,\dots,E_n)$ are \emph{equivalent} if there is $k\in\Z$ with $C_j=E_{j+k}$, where the indices run modulo $n$.
An equivalence class of closed rank-one galleries is called a \emph{loop} of galleries.

Let $\CC_1$ denote the set of all primitive closed rank-one geodesics in $\Ga\bs X$ modulo reversal of orientation.
Let $\CC_2$ denote the set of all primitive loops of galleries in $\Ga\bs X$ modulo reversal of orientation.
Let 
$$ 
Z_j(u)\df \prod_{c\in\CC_j} (1-u^{l(c)})
$$
for $j=1,2$.

\begin{proposition}\label{2.3}
For the zeta function $Z$ we have $\ds Z(u)=\frac{Z_1(u)}{Z_2(u)}$. 
Moreover, if $\Ga$ is regular, then $Z_2(u)=1$.
\end{proposition}

\prf
For any two topological spaces $X,Y$ let $[X,Y]$ be the set of homotopy classes of continuous maps from $X$ to $Y$.
Let $S^1$ be the 1-sphere and consider the natural bijection
$$
\Ga/{\rm conjugation} \ \to\ [S^1,\Ga\bs X]
$$
given by the identification $\Ga\cong\pi_1(\Ga\bs X)$.
If two closed geodesics $c_1,c_2$ are homotopic, then they are closed by conjugate elements of $\Ga$.
So they have preimages $\tilde c_1,\tilde c_2$ in $X$ which are closed by the same element $\ga$.
Hence $\tilde c_1$ and $\tilde c_2$ lie both in the apartment $\a$ where $d_\ga(x)$ is minimized.
Since $\langle\ga\rangle\bs\a$ is a cylinder, $c_1$ and $c_2$ are homotopic through closed geodesics of the same length passing through loops of galleries or intermediate rank-one geodesics.

On the other hand, each closed loop of galleries in $\Ga\bs X$ induces a homotopy between two closed geodesics of the same length: the two boundary components of the gallery.
Thus we see that the overcounting in $Z_1(u)$ is remedied by dividing by $Z_2(u)$ to result in $Z(u)$.

For the second part assume there is a closed loop $l$.
Let $(C_0,\dots,C_n)$ be a gallery in $X$ being mapped to $l$.
Then there is $\ga\in\Ga$ with $\ga C_0=C_n$ and $C_0\subset M_\ga$, where 
$$
M_\ga\=\{ x\in X : d_\ga(x)={\rm min}\}.
$$
For any $\ga$ the set $M_\ga$ is either a geodesic line or an apartment.
Since $C_0\subset M_\ga$ is our given case, it follows that $M_\ga$ is an apartment attached to a maximal split torus $A$ which contains $\ga$.
But since $\ga$ translates along a rank-one geodesic it must lie in a one-dimensional standard subtorus of $A$, which means it is not regular.
A contradiction. The claim follows.
\qed

\section{The zeta function for the $1$-skeleton}

Let $E(X)$, resp. $E(\Ga\bs X)$ denote the set of positively oriented edges in $X$, resp. $\Ga\bs X$.

Consider the vector spaces
$$
C_1(X) \= \prod_{e\in E(X)}\C e,\qquad C_1(\Ga\bs X) \= \prod_{e\in E(\Ga\bs X)}\C e.
$$
The second space is finite dimensional.
This notion actually makes sense due to Lemma \ref{orientation}.
Define a linear operator $T$ on $C_1(\Ga\bs X)$ by $Te=\sum_{e': e\to e'} e'$,
where the sum runs over all positively oriented edges $e'$ such that the endpoint of $e$ is the starting point of $e'$ and $e,e'$ lie on a rank-one geodesic, i.e., have angle $\pi$.
By the same formula, we define an operator $\tilde T$ on $C_1(X)$.
Note that $\Ga$ acts on $C_1(X)$ and that $\tilde T$ is $\Ga$-equivariant.
One has a natural identification 
$C_1(\Ga\bs X)\ \cong\ C_1(X)^\Ga$, and $T\cong \tilde T|_{C_1(X)^\Ga}$.

\begin{theorem}
We have $Z_1(u)=\det(1-uT)$. In particular, $Z_1(u)$ is a polynomial of degree equal to the number of edges of $\Ga\bs X$, or, equivalently,
$$
\deg Z_1(u)\= \frac{(q+1)N}2.
$$
where $N$ is the number of vertices in $\Ga\bs X$.
\end{theorem}

\prf
One computes
$$
\tr T^n\=\sum_e\sp{Te,e}\=\sum_{c: l(c)=n} l(c_0),
$$
where the second sum runs over all closed geodesics of length $n$ and $c_0$ is the underlying primitive of $c$.
In the next computation, we will use the letter $c$ for an arbitrary closed geodesic, $c_0$ for a primitive one, and if both occur, it will be understood that $c_0$ is the primitive underlying $c$.
We compute
\begin{eqnarray*}
Z_1(u) &=& \exp\(-\sum_{c_0}\sum_{m=1}^\infty \frac{u^{l(c_0)m}}m\)\\
&=& \exp\(-\sum_c\frac{u^{l(c)}}{l(c)}l(c_0)\)\\
&=& \exp\(-\sum_{n=1}^\infty \frac{u^n}n\sum_{c:l(c)=n}l(c_0)\)\\
&=& \exp\(-\sum_{n=1}^\infty \frac{u^n}n\tr T^n\)\\
&=& \det (1-uT).
\end{eqnarray*}
For the last line we used the fact that for a matrix $A$ we have $\exp(\tr(A))=\det(\exp(A))$.
To prove the final assertion of the Theorem it suffices to show that $T$ is invertible on $C_1(\Ga\bs X)$.
For this in turn it suffices to show that $\tilde T$ has a right-inverse on $C_1(X)$.
So let $e$ be a positively oriented edge with endpoint $[L_0]$, and let $e'$ be a positively oriented edge with start point $[L_0]$ such that $e,e'$ lie on a geodesic.
Let $[\La_2]$ be the start point of $e$ and $[\La_1]$ the end point of $e'$.
The situation is this:
$$
\begin{diagram} 
\node{[\La_2]} \arrow{e,t}{e}
	\node{[L_0]} \arrow{e,t}{e'}
		\node{[\La_1],}
\end{diagram}
$$
where $[\La_j]$ is of type $j$ for $j=1,2$.
We can choose representatives satisfying
$$
\varpi L_0\ \subset\ \La_1,\La_2\ \subset\ L_0.
$$
The condition on the types translates to the $\F_q$-vector space $\La_j/\varpi L_0$ being of dimension $j$.
The condition that $e,e'$ lie on a geodesic is equivalent to $\La_1\subset\hspace{-10pt}/\hspace{5pt}\La_2$.

For $j=1,2$ let $W_j$ be the complex vector space formally spanned by the set of all $j$-dimensional sub vector spaces of $\F_q^3$.
Let
$$
T\colon W_2\ \to\ W_1
$$
be given by
$$
T(\La_2)\=\sum_{\La_1\subset\hspace{-5pt}/\hspace{2pt}\La_2}\La_1.
$$
Define $T'\colon W_1\to W_2$ by
$$
T'(\La_1)\= \frac{-1}{q+1}\sum_{\La_2\supset\La_1}\La_2\ +\ \frac 1{q^2-q-1}\sum_{\La_2\supset\hspace{-5pt}/\hspace{2pt}\La_1}\La_2.
$$
Then
\begin{eqnarray*}
\tilde T T'(\La_1) &=& \frac{-1}{q+1}\sum_{\La_2\supset\La_1}\sum_{\La_1'\subset\hspace{-5pt}/\hspace{2pt}\La_2}\La_1'\ +\ \frac 1{q^2-q-1}\sum_{\La_2\supset\hspace{-5pt}/\hspace{2pt}\La_1}\sum_{\La_1'\subset\hspace{-5pt}/\hspace{2pt}\La_2}\La_1'\\
&=& \sum_{\La_1'} c(\La_1')\La_1',
\end{eqnarray*}
where 
\begin{eqnarray*}
c(\La_1') &=& \frac{-\#\{\La_2\supset\La_1, \La_2\supset\hspace{-11pt}/\hspace{5pt}\La_1'\}}{q+1}+\frac{\#\{\La_2\supset\hspace{-11pt}/\hspace{5pt}\La_1, \La_2\supset\hspace{-11pt}/\hspace{5pt}\La_1'\}}{q^2-q-1}\\
&=&\begin{cases} 1 & {\rm if}\ \La_1'=\La_1\\ 0 & {\rm if}\ \La_1'\ne\La_1.\end{cases}
\end{eqnarray*}
This calculation shows that the operator $T'$ on $C_1(X)$ given by
$$
T'(e)\=\frac{-1}{q+1}\sum_{\stackrel{e'\to e}{\rm non\ geodesic}}e'+\frac 1{q^2-q-1}\sum_{\stackrel{e'\to e}{\rm geodesic}} e'
$$
is a right-inverse to $\tilde T$.
The claim follows.
\qed

\subsection{A combinatorial computation}
In the following, we will write $c$ for an arbitrary closed rank-one geodesic in $\Ga\bs X$ and $c_0$ for a primitive one.
If $c$ and $c_0$ both occur, it will be understood that $c_0$ is the underlying primitive of $c$.
We compute
\begin{eqnarray*}
\frac{Z_1'}{Z_1}(u) &=& \(\log Z_1(u)\)'\\
&=& -\sum_{c_0}\sum_{n=1}^\infty l(c_0)u^{l(c_0)n-1}\\
&=& -\sum_{n=1}^\infty u^{n-1} \sum_{c: l(c)=n}l(c_0).
\end{eqnarray*}
Note that the sums run modulo reversal of orientation.

There is a natural orientation on each rank-one geodesic in $X$ given as follows.
We say that a rank-one geodesic $C$ in $X$ is \emph{positively oriented} if it runs through the vertices in the order of types: $0,1,2,0,1,2,\dots$.
The image $c_\Ga$ of $c$ in $\Ga\bs X$ is isomorphic to the image in $\langle\ga\rangle\bs X$, where $\ga\in\Ga$ is the element in $\Ga$ that closes $c$.
Since $\ga c=c$ and $\ga$ acts on $c$ by a translation it preserves the orientation of $c$ and so it does make sense to speak of positive or negative orientation for $c_\Ga$.

A \emph{line segment} in $X$ is a sequence of vertices $s=(x_0,\dots,x_k)$ such
that they are consecutive vertices on a rank-one geodesic.
A line segment in $\Ga\bs X$ is the image of one in $X$.
The \emph{length} of a line segment $s=(x_0,\dots,x_n)$ is $l(s)=n$.
On the vector spaces
$$
C_0(X)\df \bigoplus_{x\ {\rm vertex\ in\ }X} \C x,\qquad
C_0(\Ga\bs X)\df \bigoplus_{x\ {\rm vertex\ in\ }\Ga\bs X} \C x,
$$
we define an operator $A_n$ for each $n\in\N$ by
$$
A_nx\=\sum_{s:l(s)=n,\ o(s)=x} e(s),
$$
where the sum runs over all positively oriented line segments $s$ in $\Ga\bs X$ with starting point $x$ and length $n$, and $e(s)$ denote the endpoint of $s$.

\begin{lemma}
The operator $A_n$ has the trace
$$
\tr A_n\= \sum_{c: l(c)=n} l(c_0),
$$
where the sum runs over all closed rank-one geodesics in $\Ga\bs X$ modulo reversal of orientation.
\end{lemma}

\prf
Instead of summing modulo reversal of orientation one can as well sum over all positively oriented geodesics.
Recall $\tr A_n=\sum_x \sp{A_n x,x}$, where the sum runs over all vertices of $\Ga\bs X$ and the pairing $\sp ,$ is the one given by $\sp{x,y}=\delta_{x,y}$ for vertices $x,y$.
A vertex $x$ can only have a non-zero contribution $\sp{A_nx,x}$ if it lies on a close geodesic of length $n$.
The contribution of each given geodesic $c$ equals $l(c_0)$.
\qed

\subsection{The unramified Hecke algebra}
Recall that $G$ is a unimodular group, so any Haar-measure will be left- and right-invariant.
We normalize the Haar measure so that the compact open subgroup $K$ gets volume $1$.
For a subset $A$ of $G$ we write $\1_A$ for its indicator function.
Let $\CH_K$ denote the space of compactly supported functions $f\colon G\to\C$ with $f(k_1 xk_2)=f(x)$ for all $k_1,k_2\in K, x\in G$.
This is an algebra under convolution,
$$
f*g(x)\=\int_G f(y)g(y^{-1}x)\,dx.
$$
Ii is known \cite{Cartier}, that $\CH_K$ is a commutative algebra.
It has a unit element given by $\1_K$.

We will also write $KgK$ for the function $\1_{KgK}\in\CH_K$.
So a typical element of $\CH_K$ is written as
$$
f\=\sum_j c_j\, Kg_j K,\quad{\rm finite\ sum},
$$
and
$$
I(f)\= \sum_j c_j\,\vol(Kg_jK).
$$
The space $C_c(G/K)$ can be identified with $C_0(X)$ since $G/K$ can be identified with the set of vertices via $gK\mapsto gL_0$.
Likewise, $C_c(\Ga\bs G/K)$ identifies with $C_0(\Ga\bs X)$.
The Hecke algebra $\CH_K$ acts on $C_c(G/K)$ and $C_c(\Ga\bs G/K)$ via $g\mapsto g*f$, $g\in C_c(G/K)$, $f\in \CH_K$.
This will be considered as a left action as is possible since $\CH_K$ is commutative.
In $\CH_K$ we consider the elements
$$
\pi_1\= K\,\diag(1,1,\varpi)\,K,\qquad \pi_2\= K\,\diag(1,\varpi,\varpi)\, K.
$$

\begin{lemma}\label{pi_j}
For $j=1,2$,
$$
\pi_j L_0\=\sum_{\stackrel{x\ {\rm adjacent\ to\ }L_0}{x\ {\rm of\ type}\ j}}x.
$$
\end{lemma}

\prf
Clear.
\qed

\begin{proposition}\label{relations} As operators on $C_0(X)$ or $C_0(\Ga\bs X)$ respectively,
\begin{enumerate}
\item 
$A_1=\pi_1$,
\item 
$A_2=\pi_1^2-(q+1)\pi_2$,
\item 
$A_3=\pi_1^3-(2q+1)\pi_1\pi_2+(1+q+q^2)q$,
\item 
For $n\ge 3$,
$$
A_{n+1}=A_n\pi_1-qA_{n-1}\pi_2+q^3 A_{n-2}.
$$
\end{enumerate}
\end{proposition}

\prf
Part (a) follows from Lemma \ref{pi_j}.
It is clear that $\pi_1^2=A_2+c\pi_2$ for some number $c$.
From (b) in Lemma \ref{graph} it follows that $c=q+1$ which implies part (b).
The rest follows similarly.
\qed

Let $F(u)$ be the following formal powers series with values in the space $\End(C_0(\Ga\bs X))$,
$$
F(u)\=\sum_{n=1}^\infty u^{n-1} A_n.
$$
Then $\tr F(u)=\frac{Z_1'}{Z_1}(u)$.
The relations in Proposition \ref{relations} imply the following Lemma.

\begin{lemma}
We have
$$
F(u)\= H(u)\( 1-u\pi_1+u^2 q\pi_2-u^3q^3\)^{-1},
$$
where $H(u)$ is the polynomial
\begin{eqnarray*}
H(u)&=& (\pi_2-\pi_1^2) +u(\pi_1^3-\pi_1\pi_2+\pi_1^2-(q+1)\pi_2)\\ &&+u^2(\pi_1^3-(2q+1)\pi_1\pi_2+(1+q+q^2)q).
\end{eqnarray*}
\end{lemma}

\prf
This follows from Proposition \ref{relations} by a straightforward computation.
\qed

\begin{theorem}\label{3.5}
There is $m\in\N$ and a polynomial $Q(u)$ such that
$$
Z_1(u)\=\frac{\det(1-u\pi_1+u^2 q\pi_2-u^3q^3)^m}{Q(u)}.
$$
\end{theorem}

\prf
We have $\frac{Z_1'}{Z_1}(u)=\tr F(u)$, so the poles of $\frac{Z_1'}{Z_1}(u)$ must be singularities of $F(u)$, which form a subset of the set of zeros of the polynomial $\det(1-u\pi_1+u^2 q\pi_2-u^3q^3)$.
This implies the claim.
\qed

\section{The zeta function on galleries}
We now will show that the zeta function on galeries, $Z_2(u)$, also is a polynomial.
Recall that every chamber $C$ of $X$ or $\Ga\bs X$ has three vertices, one of each type $0,1,2$.
Accordingly, it has three edges of types $(0,1), (1,2)$, and $(2,0)$ respectively.
So let
$$
C_2(X)\=\prod_C \C C,\quad C_2(\Ga\bs X)\= \prod_{C\mod \Ga} \C C,
$$
where the product runs over all chambers of the buildings $X$ and $\Ga\bs X$.
On  $C_1(\Ga\bs X)$ we define a linear operator $L_1$ mapping a chamber $C$ to the sum of all chambers $C'$ such that the (1,2)-edge of $C'$ is the direct geodesic prolongation of the (0,1)-edge of $C$ as in the following picture.

$ $

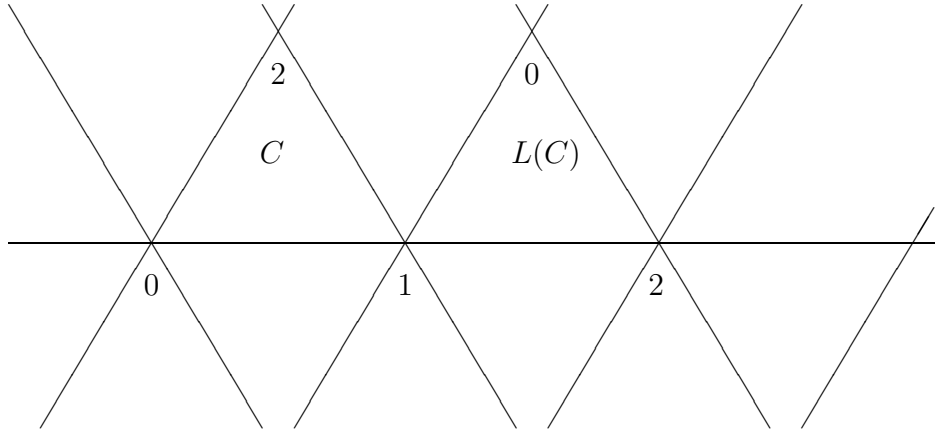
\begin{picture}(300,160)
% waagrecht
%\put(0,30){\line(1,0){350}}
\put(0,70){\line(1,0){350}}
%\put(0,110){\line(1,0){350}}
%\put(0,150){\line(1,0){350}}

% fallend
%\put(0,80){\line(3,-5){48}}
\put(0,160){\line(3,-5){96}}
%\put(48,160){\line(3,-5){96}}
\put(96,160){\line(3,-5){96}}
%\put(144,160){\line(3,-5){96}}
\put(192,160){\line(3,-5){96}}
%\put(240,160){\line(3,-5){96}}
%\put(288,160){\line(3,-5){62}}
%\put(336,160){\line(3,-5){14}}

% steigend
%\put(0,140){\line(3,5){12}}
%\put(0,60){\line(3,5){60}}
\put(12,0){\line(3,5){96}}
%\put(60,0){\line(3,5){96}}
\put(108,0){\line(3,5){96}}
%\put(156,0){\line(3,5){96}}
\put(204,0){\line(3,5){96}}
%\put(252,0){\line(3,5){96}}
\put(300,0){\line(3,5){50}}

% Bezeichnungen
\put(95,100){$C$}
\put(190,100){$L(C)$}
\put(51,50){0}
\put(147,50){1}
\put(242,50){2}
\put(99,130){2}
\put(195,130){0}

\end{picture}

$ $

\noindent
Similarly, define $L_2$ and $L_3$ by replacing $(0,1,2)$ by $(1,2,0)$ and $(2,0,1)$ respectively.
Then let $L\df L_3L_2L_1$.

\begin{proposition}\label{gall}
We have
$$
Z_2(u)\= \det(1-u^3 L).
$$
In particular, $Z_2(u)$ is a polynomial of degree at most $3$ times the number of chambers of $\Ga\bs X$.
\end{proposition}

\prf
It is easy to see that
$$
\tr L^n \= \sum_{c: l(c)=3n} \frac{l(c_0)}{3},
$$
where the sum runs over all loops of galleries in $\Ga\bs X$.
From this, the proposition follows by the same computation as before.
\qed

Finally, Theorem \ref{main} follows from Proposition \ref{2.3} , Theorem \ref{3.5} and Proposition \ref{gall}.

\noindent
{\small Mathematisches Institut\\  Auf der Morgenstelle 10\\ 72076 T\"ubingen\\ Germany\\ {\tt deitmar@uni-tuebingen.de}\\ \\
Department of Mathematics\\
Louisiana State University\\
Baton Rouge, LA 70803-4918\\ USA\\ \tt hoffman@math.lsu.edu}


\begin{thebibliography}{XXX}

\bibitem{Bass}
 \bf Bass, H.:
 \it The Ihara-Selberg zeta function of a tree lattice.
 \rm Int. J. Math. 3, No.6, 717-797 (1992).

\bibitem{Borel}
 \bf Borel, A.:
 \it Introduction aux groupes arithm\'etiques.
 \rm Hermann, Paris 1969.

\bibitem{Brown}
\bf Brown, K.:
\it Buildings. 
\rm Springer-Verlag, New York, 1989.

\bibitem{Cartier}
 \bf Cartier, P.:
 \it Representations of $\mathfrak p$-adic groups: A survey.
 \rm Automorphic forms, representations and L-functions, Proc. Symp. Pure Math. Am. Math. Soc., Corvallis/Oregon 1977, Proc. Symp. Pure Math. 33, 1, 111-155 (1979).

\bibitem{padgeom}
 \bf Deitmar, A.:
 \it Geometric zeta-functions on p-adic groups.
 \rm Math. Japon. 47, No. 1, 1-17 (1998).

\bibitem{Hash1}
 \bf Hashimoto, K.:
 \it Zeta functions of finite graphs and representations of p-adic groups.
 \rm Automorphic forms and geometry of arithmetic varieties. Adv. Stud. Pure Math. 15, 211-280 (1989).
 
\bibitem{Hash2}
 \bf Hashimoto, K.:
 \it On zeta and L-functions on finite graphs.
 \rm Int. J. Math. 1. no 4, 381-396 (1990)
 
  \bibitem{Hash3}
 \bf Hashimoto, K.:
 \it Artin type L-functions and the density theorem for prime cycles on finite graphs.
 \rm Int. J. Math. 3, no 6, 809-826 (1992).
 
\bibitem{Ihara}
\bf Ihara, Y.: 
\it Discrete subgroups of ${\rm PL}(2,\,k\sb{\wp })$.
\rm Algebraic Groups and Discontinuous Subgroups (Proc. Sympos. Pure Math., Boulder,
Colo., 1965) pp. 272--278 Amer. Math. Soc., Providence, R.I. (1966).

\bibitem{ST}
\bf Stark, H. M.; Terras, A. A.:
\it Zeta functions of finite graphs and coverings.
\rm Adv. Math. 121 (1996), no. 1, 124--165.

\end{thebibliography}
\end{document}